\documentclass[12pt]{article}

\usepackage[affil-it]{authblk}
      \usepackage{latexsym}
         \usepackage[reqno, namelimits, sumlimits]{amsmath}
         \usepackage{amssymb, amsfonts}
         \usepackage{amsthm}
 \newtheorem{theorem}{Theorem}[section]

\title{A Liouville Type Theorem for Steady-State Navier-Stokes Equations}
\author{ G Seregin
  }
  
\affil{ University of Oxford, UK and PDMI, RAS, Russia}
\date{ \today}

\begin{document}
\maketitle
\begin{abstract}
\end{abstract}
\setcounter{equation}{0}
A Liouville type theorem  is proven for the steady-state Navier-Stokes equations. It follows from the corresponding theorem on the Stokes equations with the drift. The drift is supposed to belong to a certain Morrey space.

\setcounter{equation}{0}
\section{The Main Result}
%\subsection{Normed Spaces}
The classical Liouville type theorem for the stationary Navier-Stokes  equations can be stated as follows: show that any bounded solution to the system
%$$\int\limits_{\mathbb R^3}|\nabla u|^2 dx<\infty,$$
\begin{equation}\label{nse}
u\cdot \nabla u-\Delta u=\nabla p, \qquad {\rm div }\,u=0
\end{equation}
is constant. This problem has not solved yet and even it is not clear if it has a positive answer.

Another popular problem is to show that any solution to system (\ref{nse}),    satisfying two conditions:
\begin{equation}\label{finite dissipation}
\int\limits_{\mathbb R^3}|\nabla u|^2dx<\infty.
\end{equation}
and
\begin{equation}
	\label{limit}
u(x)\to 0\quad {\rm as} \quad|x|\to \infty,
\end{equation}
is identically equal to zero. Unfortunately, it is still unknown whether the this statement is true or not.
 
 However, some attempts have been made to solve above or related problems. One of the best results in that direction can be found in \cite{Galdi-book} where it is shown that  the assumption
\begin{equation}\label{9/2}
u\in L_\frac 92(\mathbb R^3)
\end{equation}
implies $u=0$. Very recently, condition (\ref{9/2}) has been improved logarithmically in \cite{ChWo2016}.

Another set of admissible functions for solutions to (\ref{nse}), in which the Liouville type theorem is valid, has been described in  \cite{Ser2016}.   To be precise, any solution to (\ref{nse}), obeying the inclusion
\begin{equation}
	\label{L6}
	u\in L_6(\mathbb R^3)\cap BMO^{-1}(\mathbb R^3),
\end{equation}
is identically equal to zero.

For more Liouville type results, we refer the  reader to interesting  papers \cite{GilWein1978}, \cite{KNSS2009}, \cite{ChaeYoneda2013}, and \cite{Chae2014} and references there.

Our short note is inspired by paper \cite{NazUr} by Nazarov-Uraltseva
about properties of solutions to elliptic and parabolic linear equations with divergence free drift. Although their approach works for scalar equations only, similar assumptions on the drift occur in the vectorial case as well. We formulate our result as a statement  of the linear theory, considering the following steady-state Stokes system with the drift
\begin{equation}\label{stokes}
u\cdot \nabla v-\Delta v=\nabla q, \qquad {\rm div }\,v=0, \quad {\rm div}\, u=0.
\end{equation}
%Our result is follows.
\begin{theorem}
\label{main result}	
Suppose that smooth functions $u$ and $v$ satisfy (\ref{stokes}) and two additional conditions:
\begin{equation}
\label{con1}	
M:=\sup\limits_{R>0} R^{1-\frac 3q}\|u\|_{L^{q,\infty}(B(R))}<\infty
\end{equation}
with $3/2<q\leq3$ and 
\begin{equation}
\label{con2}	
N:=\sup\limits_{R>0} R^{\frac 12-\frac 3s}\|v\|_{s,B(R)}<\infty
\end{equation}
with $2\leq s\leq6$. Then $v\equiv0$ in $\mathbb R^3$.
\end{theorem}

Here, $L^{q,\infty}(\Omega)$ stands for a weak Lebesgue space, which is a particular Lorentz space $L^{q,r}(\Omega)$ and $L^{q,q}(\Omega)=L_q(\Omega)$ is a usual Lebesgue space.

It is an interesting question to understand difference between above conditions (\ref{9/2}) and (\ref{con1}), (\ref{con2}) for $u=v$. To this end, assume that there exists a divergence free field $u$ having the following bound from above
$$|u(x)|\leq \frac 1{|x'|+1}\frac 1{(|x_3|+1)^\frac 29}.$$
Then condition (\ref{con1}) holds if $q$ is not equal to 2
and condition (\ref{con2}) holds with $s=6$ while condition (\ref{9/2}) is violated.

\setcounter{equation}{0}
\section{Proof of Main Result}
\subsection{Caccioppoli Type Inequality}

Let $0<R<2$ and let a non-negative cut-off function $\varphi\in C^\infty_0(B(R))$ satisfy the following properties: $\varphi(x)=1$ in $B(r)$, $\varphi(x)=0$ out of $B(R)$, and $|\nabla \varphi(x)|\leq c/(R-r)$ for any $1\leq r<R\leq 2$. We let 
$  \overline{u}=v-[v]_{B(R)}$, where $[v]_{B(R)}$ is the mean value of $v$ over the ball of radius $R$ centred at the origin.

A given exponent $q$, satisfying conditions of Theorem \ref{main result}, see (\ref{con1}), one can find a constant $c_0(q)$ and a function $w_R$ that is smooth in $B(2)$, vanishes outside $B(R)$ and satisfies the identity  ${\rm div}\,w_R=\nabla \varphi\cdot \overline v$ and the inequality
\begin{equation}\label{Bogovskii}
\|\nabla w_R\|_{L^{2q',2}(B(R))}\leq c_0\|\nabla \varphi\cdot \overline v\|_{L^{2q',2}(B(R))}\leq \frac {c_0} {R-r}\| \overline v\|_{L^{2q',2}(B(R))}.\end{equation}
Moreover, by interpolation and Hardy-Littlewood-Sobolev inequality, we also have a bound for the right hand side of (\ref{Bogovskii}):
\begin{equation}
	\label{multiplicative}
	\| \overline v\|_{L^{2q',2}(B(R))}<c(q)\|\overline{v}\|_{2,B(R)}^{1-\frac 3{2q}}\|\nabla v\|^{\frac 3{2q}}_{2,B(R)}.                                                                                                                                                                                                                                                                                                      
\end{equation}

Now, let us test the first  equation in (\ref{stokes}) 
  with the function $\varphi \overline v-w_R$, integrate by parts in $B(R)$, and find the following identity
$$\int\limits_{B(R)}\varphi |\nabla v|^2dx=
-\int\limits_{B(R)}\nabla v :(\nabla \varphi\otimes\overline v) dx+\int\limits_{B(R)}\nabla w_R :\nabla v dx+$$
$$-\int\limits_{B(R)}(u\cdot\nabla v)\cdot\varphi\overline v dx+\int\limits_{B(R)}(u\cdot\nabla v)\cdot w_Rdx=I_1+I_2+I_3+I_4.$$

$I_1$ can be estimated easily. As a result, the below bound is valid:
$$|I_1|\leq \frac c{R-r}\|\nabla v\|_{2,B(R)}\|\overline v\|_{2,B(R)}. $$

As to $I_2$, by H\"older inequality,
we have
$$|I_2|\leq \|\nabla v\|_{2,B(R)}\|\nabla w_R\|_{2,B(R)}=\|\nabla v\|_{2,B(R)}\|\nabla w_R\|_{L^{2,2}(B(R))}\leq $$
$$\leq \|\nabla v\|_{2,B(R)}\|\nabla w_R\|_{L^{2q',2}(B(R))}\|1\|_{L^{{2q},\infty}(B(R))}\leq $$$$\leq cR^{\frac 3{2q}}\|\nabla v\|_{2,B(R)}\|\nabla w_R\|_{L^{2q',2}(B(R))}.$$
Now, taking into acount (\ref{Bogovskii}) and (\ref{multiplicative}), one can derive from the latter estimate the following:
$$|I_2|\leq c\|\nabla v\|_{2,B(R)}\frac {R^{\frac 3{2q}}}{R-r}	\|\overline{v}\|_{2,B(R)}^{1-\frac 3{2q}}\|\nabla v\|_{2,B(R)}^{\frac 3{2q}}\leq $$
$$\leq c\frac R{R-r}\|\nabla v\|^{1+\frac 3{2q}}_{2,B(R)} \Big(\frac 1{R}\|\overline {v}\|_{2,B(R)}\Big)^{1-\frac 3{2q}}.$$

Let us start evaluation of $I_3$ with integration by parts that gives
$$I_3=\frac 12\int\limits_{B(R)}|\overline{v}|^2u\cdot \nabla \varphi dx.$$
Hence, 
$$|I_3|\leq \frac c{R-r}\|u\|_{L^{q,\infty}(B(R))}\||\overline{v}|^2\|_{L^{q',1}(B(R))}\leq  $$$$\leq\frac c{R-r}\|u\|_{L^{q,\infty}(B(R))}\|\overline{v}\|^2_{L^{2q',2}(B(R))}\leq $$$$\leq \frac c{R-r}\|u\|_{L^{q,\infty}(B(R))}\|\overline{v}\|^{2(1-\frac 3{2q})}_{2,B(R)}\|\nabla v\|^{2\frac 3{2q}}_{2,B(R)}\leq $$
$$\leq c\frac R{R-r}M_0\Big(\frac 1{R^2}\|\overline{v}\|^2_{2,B(R)}\Big)^{1-\frac 3{2q}}\|\nabla v\|^{2\frac 3{2q}}_{2,B(R)}, $$
where
$$M_0=\sup\limits_{0<R<2}R^{1-\frac 3{2q}}\|u\|_{L^{q,2}(B(R))}.$$

The last term can be estimated in a similar way. Indeed, integrating by parts and applying H\"older inequality,
$$|I_4|=\Big|
\int\limits_{B(R)}(u\cdot\nabla w_R)\cdot \overline {v}dx\Big|\leq \|u\|_{L^{q,\infty}(B(R))}\||\nabla w_R||\overline{v}|\|_{L^{q',1}(B(R))}\leq$$
$$\leq  \|u\|_{L^{q,\infty}(B(R))}\|\nabla w_R\|_{L^{2q',2}(B(R))}\|\overline{v}\|_{L^{2q',2}(B(R))}\leq $$
$$\leq  \frac c{R-r}\|u\|_{L^{q,\infty}(B(R))}\||\overline{v}|^2\|_{L^{q',1}(B(R))}.
$$
The right hand side of the latter inequality  has been already estimated. Hence, we find
$$|I_4|\leq c\frac R{R-r}M_0\Big(\frac 1{R^2}\|\overline{v}\|^2_{2,B(R)}\Big)^{1-\frac 3{2q}}\|\nabla v\|^{2\frac 3{2q}}_{2,B(R)}.$$

Summarising four above estimates, we show
$$f(r)\leq c\frac R{R-r}f^\frac 12(R)\Big(\frac 1{R^2}\|\overline{v}\|^2_{2,B(R)}\Big)^\frac 12+$$$$
+c\frac R{R-r}(f(R))^{\frac 12(1+\frac 3{2q})} \Big(\frac 1{R^2}\|\overline {v}\|^2_{2,B(R)}\Big)^{\frac 12(1-\frac 3{2q})}+$$
$$+c\frac R{R-r}M_0\Big(\frac 1{R^2}\|\overline{v}\|^2_{2,B(R)}\Big)^{1-\frac 3{2q}}(f(R))^{\frac 3{2q}},$$
where 
$$f(R)=\|\nabla v\|^2_{2,B(R)}.$$
For any $1\leq R\leq 2$,
$$\frac 1{R^2}\|\overline{v}\|^2_{2,B(R)}\leq\|\widehat v\|^2_{2,B(2)}$$
 with $\widehat v=v-[v]_{B(2)}$.

Given $\varepsilon>0$, applying Young inequality, we find
$$f(r)\leq \varepsilon f(R)+c(M_0,q,\varepsilon)\|\widehat v\|^2_{2,B(2)}\Big(\frac 1{(R-r)^2}+
\frac 1{(R-r)^{\kappa_1}}+\frac 1{(R-r)^{\kappa_2}}\Big)$$
for any $1\leq R\leq 2$, where
$$\kappa_1=\frac 1{\frac 12(1-\frac 3{2q})},
\quad \kappa_2=\frac 1{1-\frac 3{2q}}.$$
As it has been shown in \cite{Giaquinta1983}, there exists a positive number $\varepsilon$ depending on $M_0$ and $q$ only such that
$$\int\limits_{B(1)}|\nabla v|^2dx\leq c(M_0,q)\int\limits_{B(2)}|v-[v]_{B(2)}|^2dx.$$

It is known  that the Navier-Stokes equations are invariant with respect to the shift and the scaling of the form
$$ v(x,t)\to \lambda v(\lambda x,\lambda^2t),\quad
q(x,t)\to \lambda^2q(\lambda x,\lambda^2t).$$
This allows us to get the required Caccioppoli type inequality
\begin{equation}
	\label{caccioppoli}
	\int\limits_{B(x_0,R)}|\nabla v|^2dx<c(M,q)\frac 1{R^2}\int\limits_{B(x_0,2R)}|v-[v]_{B(x_0,2R)}|^2dx
\end{equation}
being valid  for any $R>0$ and $x_0\in \mathbb R^3$.

\subsection{Proof of Theorem \ref{main result}}
We can put $x_0=0$ and use the following simple inequality
$$\frac 1{R^2}\int\limits_{B(2R)}|v-[v]_{B(2R)}|^2dx
\leq c\frac 1{R^2}\int\limits_{B(2R)}|v|^2dx\leq
\frac 1{R^{2(\frac 3s-\frac 12)}}\|v\|^2_{B(2R)}\leq cN^2$$
for any $R>0$. Passing $R\to\infty$, we conclude that 
$$\int\limits_{\mathbb R^3}|\nabla v|^2dx<\infty $$
The rest of the proof is the same as in \cite{Ser2016}.

\end{document}